\magnification 1200

%
%
\newdimen\FigSize	\FigSize=.9\hsize 
%
\newskip\abovefigskip	\newskip\belowfigskip
\gdef\epsfig#1;#2;{\par\vskip\abovefigskip\penalty -500
   {\everypar={}\epsfxsize=#1\noindent
    \centerline{\epsfbox{#2}}}%
    \vskip\belowfigskip}%
%
\newskip\figtitleskip
\gdef\tepsfig#1;#2;#3{\par\vskip\abovefigskip\penalty -500
   {\everypar={}\epsfxsize=#1\noindent
    \vbox
      {\centerline{\epsfbox{#2}}\vskip\figtitleskip
       \centerline{\figtitlefont#3}}}%
    \vskip\belowfigskip}%
%
\newcount\FigNr	\global\FigNr=0
\gdef\nepsfig#1;#2;#3{\global\advance\FigNr by 1
   \tepsfig#1;#2;{Figure\space\the\FigNr.\space#3}}%
%
%
%
\gdef\ipsfig#1;#2;{
   \midinsert{\everypar={}\epsfxsize=#1\noindent
	      \centerline{\epsfbox{#2}}}%
   \endinsert}%
%
\gdef\tipsfig#1;#2;#3{\midinsert
   {\everypar={}\epsfxsize=#1\noindent
    \vbox{\centerline{\epsfbox{#2}}%
          \vskip\figtitleskip
          \centerline{\figtitlefont#3}}}\endinsert}%
%
\gdef\nipsfig#1;#2;#3{\global\advance\FigNr by1%
  \tipsfig#1;#2;{Figure\space\the\FigNr.\space#3}}%
\newread\epsffilein    
\newif\ifepsffileok    
\newif\ifepsfbbfound   
\newif\ifepsfverbose   
\newdimen\epsfxsize    
\newdimen\epsfysize    
\newdimen\epsftsize    
\newdimen\epsfrsize    
\newdimen\epsftmp      
\newdimen\pspoints     
\pspoints=1bp          
\epsfxsize=0pt         
\epsfysize=0pt         
\def\epsfbox#1{\global\def\epsfllx{72}\global\def\epsflly{72}%
   \global\def\epsfurx{540}\global\def\epsfury{720}%
   \def\lbracket{[}\def\testit{#1}\ifx\testit\lbracket
   \let\next=\epsfgetlitbb\else\let\next=\epsfnormal\fi\next{#1}}%
\def\epsfgetlitbb#1#2 #3 #4 #5]#6{\epsfgrab #2 #3 #4 #5 .\\%
   \epsfsetgraph{#6}}%
\def\epsfnormal#1{\epsfgetbb{#1}\epsfsetgraph{#1}}%
\def\epsfgetbb#1{%
%
%
\openin\epsffilein=#1
\ifeof\epsffilein\errmessage{I couldn't open #1, will ignore it}\else
%
%
   {\epsffileoktrue \chardef\other=12
    \def\do##1{\catcode`##1=\other}\dospecials \catcode`\ =10
    \loop
       \read\epsffilein to \epsffileline
       \ifeof\epsffilein\epsffileokfalse\else
%
%
          \expandafter\epsfaux\epsffileline:. \\%
       \fi
   \ifepsffileok\repeat
   \ifepsfbbfound\else
    \ifepsfverbose\message{No bounding box comment in #1; using defaults}\fi\fi
   }\closein\epsffilein\fi}%
%
%
\def\epsfsetgraph#1{%
   \epsfrsize=\epsfury\pspoints
   \advance\epsfrsize by-\epsflly\pspoints
   \epsftsize=\epsfurx\pspoints
   \advance\epsftsize by-\epsfllx\pspoints
%
%
   \epsfxsize\epsfsize\epsftsize\epsfrsize
   \ifnum\epsfxsize=0 \ifnum\epsfysize=0
      \epsfxsize=\epsftsize \epsfysize=\epsfrsize
%
%
     \else\epsftmp=\epsftsize \divide\epsftmp\epsfrsize
       \epsfxsize=\epsfysize \multiply\epsfxsize\epsftmp
       \multiply\epsftmp\epsfrsize \advance\epsftsize-\epsftmp
       \epsftmp=\epsfysize
       \loop \advance\epsftsize\epsftsize \divide\epsftmp 2
       \ifnum\epsftmp>0
          \ifnum\epsftsize<\epsfrsize\else
             \advance\epsftsize-\epsfrsize \advance\epsfxsize\epsftmp \fi
       \repeat
     \fi
   \else\epsftmp=\epsfrsize \divide\epsftmp\epsftsize
     \epsfysize=\epsfxsize \multiply\epsfysize\epsftmp   
     \multiply\epsftmp\epsftsize \advance\epsfrsize-\epsftmp
     \epsftmp=\epsfxsize
     \loop \advance\epsfrsize\epsfrsize \divide\epsftmp 2
     \ifnum\epsftmp>0
        \ifnum\epsfrsize<\epsftsize\else
           \advance\epsfrsize-\epsftsize \advance\epsfysize\epsftmp \fi
     \repeat     
   \fi
%
%
   \ifepsfverbose\message{#1: width=\the\epsfxsize, height=\the\epsfysize}\fi
   \epsftmp=10\epsfxsize \divide\epsftmp\pspoints
   \vbox to\epsfysize{\vfil\hbox to\epsfxsize{%
      \includegraphics{#1}%
      \hfil}}%
\epsfxsize=0pt\epsfysize=0pt}%
%
%
{\catcode`\%=12 \global\let\epsfpercent=
%
%
\long\def\epsfaux#1#2:#3\\{\ifx#1\epsfpercent
   \def\testit{#2}\ifx\testit\epsfbblit
      \epsfgrab #3 . . . \\%
      \epsffileokfalse
      \global\epsfbbfoundtrue
   \fi\else\ifx#1\par\else\epsffileokfalse\fi\fi}%
%
%
\def\epsfgrab #1 #2 #3 #4 #5\\{%
   \global\def\epsfllx{#1}\ifx\epsfllx\empty
      \epsfgrab #2 #3 #4 #5 .\\\else
   \global\def\epsflly{#2}%
   \global\def\epsfurx{#3}\global\def\epsfury{#4}\fi}%
%
%
\def\epsfsize#1#2{\epsfxsize}%
%
%

\epsfverbosetrue			
\abovefigskip=\baselineskip		
\belowfigskip=\baselineskip		
\global\let\figtitlefont\bf		
\global\figtitleskip=.5\baselineskip	

\font\tenmsb=msbm10   
\font\sevenmsb=msbm7
\font\fivemsb=msbm5
\newfam\msbfam
\textfont\msbfam=\tenmsb
\scriptfont\msbfam=\sevenmsb
\scriptscriptfont\msbfam=\fivemsb
\def\Bbb#1{\fam\msbfam\relax#1}
\let\nd\noindent 
\def\qed{\hbox{\hskip 6pt\vrule width6pt height7pt depth1pt \hskip1pt}}
\def\natural{{\rm I\kern-.18em N}}
\def\a{\alpha}
\def\b{\beta}

\def\integer{{\rm Z\kern-.32em Z}}
\def\chix{{\raise.5ex\hbox{$\chi$}}}
\def\Z{{\Bbb Z}}
\def\real{{\rm I\kern-.2em R}}
\def\R{{\Bbb R}}
\def\complex{\kern.1em{\raise.47ex\hbox{
            $\scriptscriptstyle |$}}\kern-.40em{\rm C}}

\def\vs#1 {\vskip#1truein}
\def\hs#1 {\hskip#1truein}

\def\Month{\ifcase\number\month \relax\or January \or February \or
  March \or April \or May \or June \or July \or August \or September
  \or October \or November \or December \else \relax\fi }
\def\date{\Month \the\day, \the\year}

  \hsize=6truein        \hoffset=.25truein 
  \vsize=8.8truein      
  \pageno=1     \baselineskip=12pt
  \parskip=4 pt         \parindent=20pt
  \overfullrule=0pt     \lineskip=0pt   \lineskiplimit=0pt
  \hbadness=10000 \vbadness=10000 
\pageno=0

\footline{\ifnum\pageno=0\hss\else\hss\tenrm\folio\hss\fi}
\hbox{}
\vskip 1truein\centerline{{\bf On 2-Generator Subgroups of SO(3)}}
\vskip .5truein\centerline{by}
\centerline{Charles Radin${}^1$ and Lorenzo Sadun${}^2$} 
\footnote{}{1\ Research supported in part by NSF Grant No. DMS-9531584\hfil}
\footnote{}{2\ Research supported in part by NSF Grant No. DMS-9626698\hfil}
\vskip .2truein\centerline{Mathematics Department}
\centerline{University of Texas}
\centerline{Austin, TX\ \ 78712}
\vs.1
\centerline{radin@math.utexas.edu and sadun@math.utexas.edu}
\vs 1
\centerline{\bf Abstract}
\bigskip

We  classify all subgroups of $SO(3)$ that are generated by two
elements, each a rotation of finite order, about axes
separated by an angle that is a rational multiple of $\pi$. In all 
cases we give a presentation of the subgroup.  In most cases the
subgroup is the free product, or the amalgamated free product, 
of cyclic groups or dihedral groups.  The relations between the 
generators are all simple consequences of standard facts about
rotations by $\pi$ and $\pi/2$.  Embedded in the subgroups are explicit
free groups on 2 generators, as used in the Banach-Tarski paradox.

\vs.7
\centerline{October 1997}
\vs.2
\centerline{Subject Classification:\ \ 51F25, 52C22}

\vfill\eject 
\nd
{\bf \S 1. Introduction}

This paper concerns 2-generator subgroups of $SO(3)$. Let $A$ and $B$
be rotations of finite order of Euclidean 3-space, about axes that are
themselves separated by an angle which is a rational multiple of $\pi$.
We call these groups ``generalized dihedral'' groups.
We are interested in
the algebraic relations between $A$ and $B$. The special cases where
the axes are orthogonal were treated in [CR, RS1]; here we solve the
general case. Most such groups are infinite and dense. The only exceptions are:
if one generator has order 1, the group is cyclic; if one generator
has order 2 and the axes are orthogonal, the group is dihedral; and if
both generators have order 4 and the axes are orthogonal, the group is
the symmetries of the cube. (The symmetries of each other Platonic
solid can also be generated by a pair of finite order rotations, but
only with axes separated by an irrational angle.)

These 2-generator subgroups appear in the theory of tilings of 
Euclidean $\R^3$. Tilings have been constructed [CR, RS2]
in which each tile appears in an infinite number of different orientations
in space.  The set of relative orientations, {\it i.e.}, the set of rotations
that make one tile parallel to another, forms a group, precisely of the
sort described here.

For each 2-generator subgroup considered here,
we give an explicit presentation.  In many cases we also
express the group as a free or amalgamated free product. In particular
we show that all relations between such rotations can be understood in
terms of a short list of simple relations involving rotations by
$\pi/2$:

\item {1.} If $A$ is a rotation by $\pi$ about one axis 
and $B$ is any rotation about an orthogonal axis
then $ABA=B^{-1}$.  Equivalently, $ABAB$ is the identity.

\item {2.} If $A$ and $B$ are rotations by $\pi/2$ about orthogonal axes, 
then $ABABAB$ is the identity.  Equivalently, $ABA=BAB$. 

\item {3.} If $A$ and $B$ are rotations by $\pi$ about axes that are 
separated by an angle $\theta$, then $BA$ is a rotation, by $2\theta$,
about an axis perpendicular to the axes of $A$ and $B$. 

We need some further notation. See Figure 1.  Let $\ell$ be the line
in the $x$-$y$ plane, through the origin and making an angle of $2 \pi
n/m$ with the $x$-axis, where $n$ and $m$ are relatively prime
positive integers and $m>2$.  Let $A=R_x^{2\pi/p}$ be rotation by $2
\pi/p$ about the $x$-axis, where $p$ is a positive integer.  Let $B =
R_\ell^{2 \pi/q}$ be rotation by $2 \pi/q$ about the line $\ell$.  $A$
and $B$ will be our primary objects of study.  We similarly define the
rotations $C=R_z^{\pi/q}$, $S=R_y^{\pi/2}$ and $U=R_z^{2\pi/m}$.
$T$ and $\hat T$ will denote various rotations
about the $x$ axis.

Let $G_{n/m}(p,q)$ be the
subgroup of $SO(3)$ generated by $A$ and $B$.  An
important special case is where the axes of $A$ and $B$ are
orthogonal: that is, ${n/m}=1/4$.  In that case we omit the subscript
and write $\hat G(p,q)
\equiv G_{1/4}(p,q)$.  Among the groups $\hat G(p,q)$, a further special
case is where $q=4$, which will play a pivotal role. 

These groups form a hierarchy.  The structure of $\hat G(p,4)$ is
considerably simpler than that of a general $\hat G(p,q)$, while the
structure of $\hat G(p,q)$ is much simpler than that of a general
$G_{n/m}(p,q)$.  However, $G_{n/m}(p,q)$ actually embeds in the
``simpler'' $\hat G(pq,m)$, while $\hat G(p,q)$ embeds in the ``simpler''
$\hat G(pq,4)$:
$$ G_{n/m}(p,q) \subset \hat G(p',q') \subset \hat G(s,4) = \hbox{ not too
complicated}. \eqno(1) $$
We can therefore use the known properties of a simpler group to deduce
the structure of the more complicated group. The cases $\hat G(3,3)$ and 
$\hat G(3,4)$ were worked out in [CR]; the structures of $\hat G(p,4)$
and $\hat G(p,q)$ were derived in [RS1]; the structure of $G_{n/m}(p,q)$ is
new.

To see the embeddings of (1) it is useful to define $G(p,q)$ as the
subgroup (conjugate to $\hat G(p,q)$) generated by rotations about the
$x$ and $z$-axes, rather than the $x$ and $y$-axes. That is, $G(p,q)$
is generated by $A$ and $C$, rather than $A$ and $B$, while $\hat
G(m,4)$ is generated by $T=R_x^{2\pi/m}$ and $S=R_y^{\pi/2}$.  Note
that, if $s$ is a multiple of $p$ and also a multiple of $q$, then
$A=T^{s/p}$, while $C=S^{-1} T^{s/q} S$, so that $G(p,q) \subset
\hat G(s,4)$. Next we show that $G_{n/m}(p,q) \subset G(r,m)$, where $r$
is any common multiple of $p$ and $q$.  The generators of
$G_{n/m}(p,q)$ are $A=R_x^{2\pi/p}$ and $B=R_\ell^{2\pi/q}$, while the
generators of $G(r,m)$ are $\hat T=R_x^{2\pi/r}$ and $U=R_z^{2\pi/m}$.  
Since $r$ is a multiple of $p$ and $q$, $A=\hat T^{r/p}$ and
$B=U^{-n}\hat T^{r/q}U^n$.

Although there are quite a few cases, depending on the values of $p$,
$q$, $m$ and $n$, the results of this paper, and of [RS1], are all of
the following form:

\nd {\bf ``Theorem''} {\it Let $G=\hat G(m,4)$, $G(p,q)$ or$\,$
$G_{n/m}(p,q)$. There is a short list of simple relations among the
generators of $G$, all of which involve rotations by multiples of $\pi/2$,
and all of which are simple consequences of 
$$R_x^\pi R_y^\theta R_x^\pi = R_y^{-\theta}, \eqno(2) $$ 
$$R_y^{\pi/2}R_x^{\pi/2}R_y^{\pi/2}= 
R_x^{\pi/2}R_y^{\pi/2}R_x^{\pi/2}, \eqno (3)$$ 
$$R_\ell^\pi R_x^{\pi}=R_z^{4 \pi n/m}. \eqno(4)$$

\nd Any other relation may be derived from these simple
relations.}

For example, powers of the generators of $G_{n/m}(p,q)$, with $p$ and
$q$ both odd, are never rotations by multiples of $\pi/2$, unless they
are rotations by $2\pi$. So we expect (and prove) that, if $p$ and $q$
are odd, $G_{n/m}(p,q)$ is the free product of $\Z_p$ and $\Z_q$,
for any $n$ and $m$.

More precise formulations of this ``theorem'' may be found below.  The
cases $G=\hat G(m,4)$ and $G(p,q)$ are discussed in \S 2 (Theorems
1 and 2).  The cases $G_{n/m}(p,q)$ fall into 3 classes, which are
the subjects of Theorems 3--5 in \S 3.  Theorems 3--5 depend on Theorem 2, 
which in turn depends on Theorem 1, which in turn depends on the
following result:

\nd {\bf Foundation Lemma [RS1]:} {\it Let $T=R_x^{2 \pi/m}$, let $S =
R_y^{\pi/2}$, and consider the word 
$$ \chix = S^W S^{a_1} T^{b_1} S^{a_2} T^{b_2} \cdots S^{a_n} T^{b_n}
S^E, \eqno(5) $$
where $W$, $E$, $a_i$, $b_i$ and $n$ are integers, $n>0$, none of the
$a_i$'s are divisible by 2, and none of the $b_i$'s are divisible by
$m/4$. Then $\chix$ is not equal to the identity.}

The following two extensions of the Foundation Lemma will also be used
repeatedly:

\nd {\bf Extension 1:} {\it  Let $T=R_x^{2 \pi/m}$, let $S =
R_y^{\pi/2}$, and let $\chix$ be as in eq~(5). If 
$W$, $E$, $a_i$, $b_i$ and $n$ are integers, $n>0$, none of the
$a_i$'s are divisible by 2, and no two consecutive 
$b_i$'s are divisible by $m/4$, then $\chix$ is not equal to the identity.}

\nd Proof: If $b_i=m/4$, $a_{i-1}=a_i=1$ and 
$b_{i-1}$ and $b_{i+1}$ are not divisible
by $m/4$, then we can shorten our word, removing the $T^{m/4}$ term 
using the identity (3):
$$ T^{b_{i-1}} S T^{m/4} S T^{b_{i+1}} = 
 T^{b_{i-1}} T^{m/4} S T^{m/4} T^{b_{i+1}} = 
T^{b_{i-1} + m/4} S T^{b_{i+1}+m/4}.
\eqno(6) $$  
Similar manipulations apply if $a_{i-1}$, $b_i$ or $a_i$ is negative.
Note that, since $b_{i\pm 1}$ is not divisible by $m/4$, neither is
$b_{i \pm 1} + m/4$. We repeat this procedure until all the factors of
$T^{\pm m/4}$ are removed, except possibly $T^{b_1}$ and $T^{b_n}$.
If $b_1$ and $b_n$ are not multiples of $m/4$, the Foundation Lemma
applies, and we are done.  If $b_1$ or $b_n$ equals $\pm m/4$ we apply
the proof of the Foundation Lemma to the terms in the middle (e.g.,
$S^{a_2} T^{b_2} \cdots S^{a_n}$ if $b_1= b_n =m/4$) and see that this
central piece has four non-integral matrix elements.  Since $S$ and
$T^{m/4}$ are, up to sign, permutation matrices, $\chix$ has four
non-integral matrix elements and cannot equal the identity. \qed

\nd {\bf Extension 2:} {\it Let $A$ and $C$ be rotations about orthogonal
axes by $2 \pi/p$ and $2 \pi/q$, respectively, 
and consider a word of the form 
$$ A^{a_1} C^{c_1} \cdots A^{a_n} C^{c_n}, \eqno(7) $$
where $a_1$ and $c_n$ may be zero, but otherwise none of the $a_i$'s
are multiples of $p/2$ and none of the $c_j$'s are multiples of $q/2$.
If no two consecutive terms represent rotations by $\pm \pi/2$ (e.g.\
$A^{p/4} C^{-q/4}$ or $C^{q/4}A^{p/4}$), and if the word is not empty,
then the word is not equal to the identity.}

\nd Proof: Choose axes such that $A=R_x^{2\pi/p}$ and $C=R_z^{2\pi/q}$, and 
let $S=R_y^{\pi/2}$ and $T=R_x^{2 \pi/pq}$, so that $A=T^q$
and $B = S^{-1} T^p S$. Rewrite the word in terms of $S$ and $T$ and apply 
Extension 1. \qed

To prove Theorems 1--5, we employ two strategies repeatedly.  We
always begin by noting certain simple relations among the generators
of $G$.  The first strategy is to use these relations to convert an
arbitrary word in the generators either into the identity or into a
word which, by the Foundation Lemma and its Extensions, is known not
to equal the identity.  When this strategy succeeds it implies that
there cannot be any relations that are independent of the stated
ones. A presentation of the group then follows.  A second strategy is
to embed the group $G$ in a seemingly larger group $H$, and then use
the simple relations to show that the generators of $H$ are actually
in $G$, and thus that $G=H$.  All cases yield to a combination of
these strategies.  In a small number of cases (Theorem 5), 
both strategies are needed.  First we choose
an appropriate $H$ and show that $G=H$, and then we apply the first
argument to $H$ to find a presentation.

\bigskip

\nd {\bf \S 2. The structures of $\hat G(m,4)$ and $G(p,q)$.}

\bigskip

Before examining the structures of $G(p,q)$ and $\hat G(m,4)$ in general, we
consider again the few cases where the group is finite.  $G(p,1)$ is the
cyclic group of order $p$, and has the presentation
$$ 
G(p,1) = \Z_p =\, <\!\a: \a^p\!>. \eqno(8) 
$$
$G(p,2)$ is the dihedral group of order $2p$, and has the presentation
$$ 
G(p,2) = D_p =\, <\!\a, \b: \a^p, \b^2, (\a\b)^2\!>. \eqno(9)
$$
The relation $(\a\b)^2=1$, which is equivalent to eq.~(2),
will be used repeatedly.  

The last finite group is $\hat G(4,4)=G(4,4)$, 
which is
the 24-element group of rotational symmetries of the cube.  It has the
presentation 
$$ 
\hat G(4,4) =\, <\!\a, \b: \a^4, \b^4, (\a\b^2)^2, (\a^2\b)^2, (\a\b)^3\!>.
\eqno(10)
$$
A new ingredient is the last relator in (10), which is a consequence of 
eq.~(3).  This, too, will be used frequently.

\nd {\bf Theorem 1:} {\it The group $\hat G(m,4)$ has the following
presentation: 

\item {1.} If $m$ is odd, then
$$
\hat G(m,4) =\, <\! \a, \b: \a^m, \b^4, (\a\b^2)^2\!> \,= D_m *_{\Z_2} \Z_4,
\eqno(11) $$

\item {2.} If $m$ is twice an odd number, then
$$
\hat G(m,4) =\, <\! \a, \b: \a^m, \b^4, (\a\b^2)^2, (\a^{m/2}\b)^2\!>\, = 
D_m *_{D_2} D_4,
\eqno(12) $$

\item {3.} If $m$ is divisible by 4, then
$$\hs-.02
\hat G(m,4) =\, <\! \a, \b: \a^m, \b^4, (\a\b^2)^2, (\a^{m/2}\b)^2, 
(\a^{m/4}\b)^3\!>\, = D_m *_{D_4} \hat G(4,4).\eqno(13) $$}

\nd Proof: Let $\hat G(m,4)$ be generated by $S=R_y^{\pi/2}$ and
$T=R_x^{2\pi/m}$. In each case, the relators listed (with $\a$
corresponding to $T$ and $\b$ to $S$), are just the previously
discussed relators either involving rotations by $\pi$ or involving
pairs of rotations by $\pi/2$.  Our task is to show there are no
additional relations.  We do this by taking an arbitrary word in $S$
and $T$ and, using the known relations, either reducing the word to
the form (5) or to the identity.  If the reduction is to the form (5),
then, by the Foundation Lemma (or Extension 1), 
the word cannot be a relator.  If the
reduction is to the identity, the word {\it is} a relator, but is
not independent of the listed relators.

We begin with case 1.  An arbitrary word in $S$ and $T$ is of the form
$$ 
S^W S^{a_1} T^{b_1} S^{a_2} T^{b_2} \cdots S^{a_n} T^{b_n}
S^E, \eqno(14)
$$
with $W$, $a_i$, $b_i$ and $E$ arbitrary.  Using the relators $S^4$
and $T^m$, we can require that the $b_i$'s be between $1$ and $m-1$,
inclusive, and that each $a_i$ be 1, 2 or 3. Since $m$ is odd,
$T^{b_i}$ can never be a rotation by a multiple of $\pi/2$. If any of
the $a_i$'s equal 2, we can shorten the word using $(TS^2)^2=1$, or
equivalently $T^bS^2 = S^2 T^{-b}$, and hence $T^b S^2 T^{b'} = S^2
T^{b'-b}$. In this way, we shorten the word until either $n=0$ or
until all the $a_i$'s are odd, in which case we have the form (5).

In case 2, $T^{b_i}$ may be a rotation by $\pi$, but not by $\pm
\pi/2$. In that case, in addition to removing $S^2$ factors, we
remove $T^{m/2}$ factors, using 
$T^{m/2}S = S^{-1} T^{m/2}$ and hence $T^{m/2} S^a = S^{-a} T^{m/2}$.  
Once the $T^{m/2}$s and $S^2$s are removed, we have the form (5). 

In case 3, even after removing the $T^{m/2}$s and $S^2$s, we may have
some $T^{\pm m/4}$s in our word.  If two consecutive factors appear 
we shorten the word using $(T^{m/4}S)^3=1$, and hence
$S T^{m/4} S T^{m/4} = T^{-m/4} S^{-1}$, with similar identities applying
if the exponents are not all positive.  
Thus the word may continue to be shortened until either Extension 1
may be applied, or the word is so short it can be checked by hand.

This gives the presentations (11--13).  As for the amalgamated free
products, in case 1 the $D_m$ subgroup is generated by $\a$ and
$\b^2$, $\Z_4$ is generated by $\b$, and their intersection is $\Z_2 =
\{ 1, \b^2\}$.  In case 2 the $D_m$ subgroup is generated by $\a$ and
$\b^2$, the $D_4$ subgroup by $\a^{m/2}$ and $\b$, and the common
$D_2$ by $\a^{m/2}$ and $\b^2$.  In case 3, $D_m$ is generated by $\a$
and $\b^2$, $\hat G(4,4)$ by $\a^{m/4}$ and $\b$, and the common $D_4$
subgroup by $\a^{m/4}$ and $\b^2$. \qed

\nd {\bf Theorem 2:} {\it The group $G(p,q)$ has the following
structure: 

\item {1.} If $p$ and $q$ are both odd, then
$$
G(p,q) =\, <\! \a, \b: \a^p, \b^q\!>\, = \Z_p * \Z_q,
\eqno(15) $$

\item {2.} If $p$ is even and $q$ is odd, then
$$
G(p,q) =\, <\! \a, \b: \a^p, \b^q, (\a^{p/2}\b)^2\!>\, = 
\Z_p *_{\Z_2} D_q,
\eqno(16) $$

\item {3.} If $p$ and $q$ are both even, but $q$ is not divisible by 
4, then
$$
G(p,q) = \,<\! \a, \b: \a^p, \b^q, (\a^{p/2}\b)^2,  (\a\b^{q/2})^2\!>\, 
=  D_p *_{D_2} D_q.
\eqno(17) $$ 

\item {4.} If $p$ and $q$ are both divisible by 4, then
$G(p,q) = \hat G([p,q],4)$.}

\nd Proof: In cases 1--3, the relations listed are known consequences
of rotations by $\pi$. What remains is to show that these
are the only relations.  

Let $A=R_x^{2\pi/p}$ and $C=R_z^{2\pi/q}$ be
the generators of $G(p,q)$.  Take an arbitrary word in $A$ and $C$, and
use the given relations as above to put it in the form
$$ A^{a_1} C^{c_1} \cdots A^{a_n} C^{c_n}, \eqno(18)$$ 
where $a_1$ and $c_n$ may be zero, but none of the other $a_i$'s are
multiples of $p/2$, and none of the other $c_j$'s are multiples of
$q/2$.  (In case 1 there is nothing to do, in case 2 we are
eliminating factors of $A^{p/2}$, and in case 3 we are eliminating
factors of $A^{p/2}$ and $C^{q/2}$).  Since $q$ is not divisible by 4, 
none of the remaining $C^{c_j}$ terms are rotations by multiples of 
$\pi/2$.  Thus, by Extension 2, our word is not the identity.

All that remains for cases 1--3 are the amalgamated free products.  In
case 2, $D_q$ is generated by $A^{p/2}$ and $C$, while $\Z_p$ is
generated by $A$.  In case 3, $D_p$ is generated by $A$ and $C^{q/2}$, 
$D_q$ is generated by $A^{p/2}$ and $C$, and their intersection,
$D_2$, is generated by $A^{p/2}$ and $C^{q/2}$. 

In case 4, we apply the second strategy.  
Consider $S$ and $\hat T = R_x^{2 \pi /[p,q]}$, the
generators of $\hat G([p,q],4)$. Since
$A$ and $C$ can be constructed from $S$ and $\hat T$, we must have
$G(p,q) \subseteq \hat G([p,q],4)$.  We will show that $S,\hat T 
\in G(p,q)$, and thus $G(p,q)=\hat G([p,q],4)$.

Since $p$ is divisible by 4, $R_x^{\pi/2} = A^{p/4} \in G(p,q)$.
Since $q$ is divisible by 4, $R_z^{\pi/2} = C^{q/4} \in G(p,q)$.
Thus $S = R_y^{\pi/2} = R_x^{-\pi/2} R_z^{\pi/2} R_x^{\pi/2} \in
G(p,q)$. Now $S^{-1} C S = R_x^{2 \pi/q} \in G(p,q)$.  Since $R_x^{2
\pi/p}$ and $R_x^{2 \pi/q}$ are in the group, so is $\hat T$.  \qed

\bigskip

\nd {\bf \S 3. The structure of $G_{n/m}(p,q)$.}

\bigskip

In studying $G(p,q)$ (Theorem 2) we saw that one of two things
happened: either the only relations among the generators $A$ and $C$
were the obvious ones, or $G(p,q)$ was equal to the previously
understood group $\hat G([p,q],4)$. In analyzing $G_{n/m}(p,q)$, we shall
have 3 possibilities.  Sometimes there are no non-trivial relations
between the generators; these cases are treated in Theorem 3.
Sometimes $G_{n/m}(p,q)= G(m,[p,q])$; these cases are treated in
Theorem 4.  In a few exceptional cases (Theorem 5), $G_{n/m}(p,q)$ is
neither $\Z_p * \Z_q$ nor $G(m,[p,q])$; these cases are in some sense
intermediate.

Our notation is as in Figure 1.  As always, $A=R_x^{2 \pi/p}$,
$S=R_y^{\pi/2}$, $\ell$ is the line in the $x$-$y$ plane, through
the origin, making an angle of $2 \pi n/m$ with the $x$ axis, and 
$B=R_\ell^{2\pi/q}$. We now fix
$T=R_x^{2 \pi/pq}$, $\hat T = R_x^{2\pi/[p,q]}$
and $U=R_z^{2 \pi /m}$. $A$ and $B$ generate $G_{n/m}(p,q)$, while $U$
and $T$ (or $\hat T$) generate $G(pq,m)$ (or $G([p,q],m)$).  Note that
$A=T^q$ and $B=U^{-n} T^p U^n$, so $G_{n/m}(p,q) \subseteq G(m,pq)$.
Similarly, $G_{n/m}(p,q) \subseteq G([p,q],m)$.  For any integer $N$,
let $\rho(N)$ be the number of powers of 2 that divide $N$.  For
example $\rho(3)=0$ and $\rho(12)=2$.

Without loss of generality we can assume that $m$ is either odd or a
multiple of 4.  For if $m$ is an odd multiple of 2, then we can
replace the axis $\ell$ by $-\ell$, which makes an angle of $2
\pi([(m/2)+n]/m)$ with the $x$ axis.  Since $m/2$ and $n$ are odd,
$m/2+n$ is even, and we can replace $n$ and $m$ by $n'=(m/2+n)/2$ and
$m'=m/2$, where $m'$ is odd.

\nd {\bf Theorem 3:} {\it If $p$ and $q$ are odd, or if $p$ is even, 
$q$ is odd and $m \ne 4$, then
$$ 
G_{n/m}(p,q) = \,<\!\a,\b: \a^p, \b^q\!>\, = \Z_p * \Z_q. \eqno(19)
$$
If $p$ is even, $q$ is odd and $m=4$, then 
$$
G_{n/m}(p,q) = \,<\! \a, \b: \a^p, \b^q, (\a^{p/2}\b)^2\!>\, = 
\Z_p *_{\Z_2} D_q.
\eqno(20) $$
}

\nd Proof: Take an arbitrary nontrivial word in $A$ and $B$,
$$ 
A^{a_1} B^{b_1} \cdots A^{a_n} B^{b_n}, \eqno(21) 
$$
and rewrite it in terms of $U$ and $T$:
$$ T^{q a_1} U^{-n} T^{pb_1} U^n \cdots T^{pb_n} U^n. \eqno(22) $$
If $p$ and $q$ are odd then no power of $T$ can be a rotation by
$\pi/2$ or $\pi$.  Since $m>2$, and since $n$ and $m$ are relatively
prime, $U^{\pm n}$ is not a rotation by $\pi$ (although it may be a
rotation by $\pi/2$ if $m=4$). By Extension 2 the word is not
equal to the identity, so no relations between $A$ and $B$ exist.

If $p$ is even, $q$ is odd, and $m \ne 4$, then $U^n$ is not a
rotation by $\pi/2$ or $\pi$.  However, some of the $T^{q a_i}$ terms
may be rotations by $\pi$.  We remove these using the fact
that 
$$ T^{pb_{i-1}} U^n T^{pq/2} U^{-n} T^{pb_i} = T^{pb_{i-1}} U^{2n}
T^{pb_i + pq/2}. \eqno(23)
$$
Note that $T^{pb_i + pq/2}$ is not a rotation by a multiple of
$\pi/2$.  There may remain some $T^{q a_i}$ terms that are rotations
by $\pi/2$, and $U^{2n}$ might be a rotation by $\pi/2$ (if $m=8$),
but these cannot occur consecutively.  Every $T^{q a_i}$ is flanked by
$U^{\pm n}$s, while every $U^{2n}$ is flanked by $T^{p b_{i-1}}$ and
$T^{pb_i+pq/2}$.  Thus by Extension 2 the word is not equal to the
identity.  

The case where $p$ is even, $q$ is odd and $m=4$ is just part 2 of
Theorem 2. \qed

What remain are the cases where both $p$ and $q$ are even.  In
these cases we have a new relation.  Rotation by $\pi$ about the $x$
axis, followed by rotation by $\pi$ about the $\ell$ axis, equals
rotation by $4 \pi n/m$ about the $z$ axis:  
$$B^{q/2} A^{p/2} = R_\ell^\pi R_x^{\pi}=R_z^{4\pi n/m}= U^{2n}. \eqno(24)$$
To see
this, note that $B^{q/2} A^{p/2}$ fixes the $z$ axis while reflecting
the $x$ axis across the $\ell$ axis.  The remainder of this paper 
tracks the effect of this relation.
\bigskip
\nd {\bf Theorem 4:} {\it If

\item {1.} $p$ is even, $q$ is even and $m$ is odd, or

\item {2.} $\rho(p)\ge 2$, $\rho(q)\ge 2$ and $\rho(m)=2$, or

\item {3.} $q$ is even and $\rho(p) \ge \rho(m) \ge 3$, 

then $G_{n/m}(p,q)=G([p,q],m)$}

\nd Proof: We must show that $U$ and $\hat T$ are in $G_{n/m}(p,q)$.
First we show it suffices to prove that $U^n \in G_{n/m}(p,q)$. Since
$n$ and $m$ are relatively prime, $U$ is a power of $U^n$.  If
$U^n \in G_{n/m}(p,q)$ then $R_x^{2\pi/q}= U^{n} B U^{-n} \in
G_{n/m}(p,q)$, and so $\hat T \in G_{n/m}(p,q)$.  Note that, by eq.~(24),
we know that $U^{2n}$, and therefore $U^2$, is in $G_{n/m}(p,q)$.

If $m$ is odd, then $U^n = (U^{2n})^{(m+1)/2}$, and case 1 is done.

In case 2, $\rho(m)=2$, so $R_z^{\pi/2}$ is an odd power of $U$, and hence an
odd power of $U^n$. Let $s$ be an odd integer such that $R_z^{\pi/2} =
U^{ns}$.  $U^{n(s-1)}$ is then an even power of $U^n$, and hence an
element of $G_{n/m}(p,q)$.  Thus $R_y^{2 \pi/q} = R_z^{-\pi/2} T^p
R_z^{\pi/2} = U^{-n(s-1)}B U^{n(s-1)} \in G_{n/m}(p,q)$. Since $4$
divides $q$, a power of this is $R_y^{\pi/2}=S$.  Since $4$ divides
$p$, $R_x^{\pi/2} = A^{p/4} \in G_{n/m}(p,q)$.  Thus $R_z^{\pi/2} =
S^{-1} R_x^{\pi/2} S \in G_{n/m}(p,q)$.  So $U^n = R_z^{\pi/2}
U^{-n(s-1)} \in G_{n/m}(p,q)$, which completes case 2.

In case 3, since $\rho(m) \ge 3$, $R_z^{\pi/2}$ is an even power of
$U^n$ and hence is in $G_{n/m}(p,q)$.  Since 4 divides $p$, $R_x^{\pi/2} \in
G_{n/m}(p,q)$. From these we see that $S\in G_{n/m}(p,q)$. Conjugating
$A$ by $S$ we obtain $R_z^{2 \pi/p} \in G_{n/m}(p,q)$.  Since $\rho(p)
\ge\rho(m)$, $U^b = [R_z^{2\pi/p}]^{a} \in G_{n/m}(p,q)$, 
where $a= p/2^{\rho(p)-\rho(m)}$ and $b=m/2^{\rho(m)}$.  Since $n$ and
$b$ are both odd, $n-b$ is even, so $U^{n-b} = (U^2)^{(n-b)/2} \in
G_{n/m}(p,q)$ and $U^n = U^b U^{n-b} \in G_{n/m}(p,q)$. \qed

All that is left are a handful of special cases, all of which have $p$
and $q$ even and $m$ divisible by 4. Without loss of generality we
assume $\rho(p)\ge \rho(q)$.

\nd {\bf Theorem 5:} {\it

\item {1.} If $\rho(m)>\rho(p)=\rho(q)=1$, then
$$
G_{n/m}(p,q) = \,<\! \a, \b, \gamma: \a^p, \b^q, \gamma^2, (\a\gamma)^2,
(\b\gamma)^2, (\a^{p/2}\b^{q/2})^{m/4}\gamma>.
\eqno(25) 
$$
This is a quotient of $D_p *_{\Z_2} D_q$ by the last relation.

\item {2.} If $\rho(m)>\rho(p) >1 $ and $\rho(m)>\rho(q)>1$, then
$$
G_{n/m}(p,q) = G_{1/2r}(4,4) = G(r,4) *_{D_{r}} G(r,4),
\eqno(26)
$$
where $r = [p,q,m/2]$ is the least common multiple of $p$, $q$ and
$m/2$, the first $G(r,4)$ subgroup is generated by
$R_z^{2\pi/r}$ and $R_x^{\pi/2}$, the second is conjugate to the first
by $R_z^{\pi/r}$, and the common $D_r$ subgroup is generated by
$R_z^{2\pi/r}$ and $R_x^\pi$.

\item {3.} If $\rho(m) > \rho(p) > \rho(q) = 1$, then 
$$
G_{n/m}(p,q) = G_{1/2s} (4,q) = G(s,4) *_{D_2} D_q,
\eqno(27)
$$
where $s=[p,m/2]$, $G(s,4)$ is generated by $R_z^{2 \pi/s}$ and
$R_x^{\pi/2}$, $D_q$ is generated by $B$ and $R_z^{\pi}$, and the common
$D_2$ subgroup is generated by $B^{q/2}$ and $R_z^\pi$.

\item {4.} If $\rho(p) > \rho(q) = 1$ and $\rho(m)=2$, then
$$
G_{n/m}(p,q) = G_{1/2t}(p,2) = D_p *_{D_2} D_{t}, \eqno(28)
$$
where $t=[q,m/2]$, $D_p$ is generated by $A$ and $R_z^{\pi}$, $D_t$ is
generated by $R_z^{2 \pi/t}$ and $B^{q/2}$, and the 
common $D_2$ is generated by $R_z^{\pi}$ and $A^{p/2}=R_z^{2 \pi/t}B^{q/2}$.}

\nd Proof:  We begin with case 1.  Taking $\a = A$, $\b = B$ and $\gamma
= R_z^\pi = (A^{p/2} B^{q/2})^{m/4}$ (from (24), since $n$ is odd), 
it is clear that $\a, \b$ and $\gamma$
generate the group (indeed, $\a$ and $\b$ do), and that all the 
stated relations hold.  As usual, we must show that there are no additional
relations.

Take an arbitrary word $w$ in $\a$, $\b$ and $\gamma$ and shorten it
as follows.  Make all exponents of $\a$ less than or equal to $p/2$ in
magnitude, and all exponents of $\b$ less than or equal to $q/2$ in
magnitude.  Use the relators $(\a\gamma)^2$ and $(\b\gamma)^2$ to move
all the $\gamma$s to the left (i.e., replace $\alpha\gamma$
by $\gamma\alpha^{-1}$ etc.), after which there should be either one
or zero $\gamma$s.  If there are any sequences $\cdots \a^{p/2}
\b^{q/2} \cdots$ with more than $m/4$ pairs, replace it with $\gamma$
times a shorter sequence, using the relator
$(\a^{p/2}\b^{q/2})^{m/4}\gamma$.  Repeat these procedures as
necessary until there is at most one $\gamma$ (at the left), all
rotations are by $\pi$ or less, and there are at most $m/4$ pairs of
rotations by $\pi$ in a row.  We claim that this word, if not
explicitly the identity, is not equal to the identity.

So consider $w = \gamma^{c}
\a^{a_1} \b^{b_1} \cdots \a^{a_n} \b^{b_n}$, with $c=0$ or $1$ and the
aforementioned restrictions on the $a_i$s and $b_i$s.
Now rewrite $\a$, $\b$ and $\gamma$ in terms of $U$ and $T$, getting
$$ w = U^{cm/2} \prod_i (T^{qa_i} U^{-n} T^{pb_i} U^n ). \eqno(29) $$
Since $p$ and $q$ are even, some of the powers of $T$ may be rotations
by $\pi$, but, since neither $p$ nor $q$ is a multiple of 4, none are
rotations by $\pm \pi/2$.  We eliminate the $T^{pq/2}$ terms, pushing
them to the left (in the same spirit as (23)). If $k$ such pairs are
adjacent, then this process changes the $k+1$ flanking $U^{\pm n}$s
into a single $U^{\pm(k+1)n}$. However, since $k$ is never bigger than
$m/4$, $U^{\pm(k+1)}$ is never a rotation by $\pi$, and since $n$ is
relatively prime to $m$ neither is $U^{\pm(k+1)n}$. The remaining
powers of $T$ are not rotation by multiples of $\pi/2$.  Thus, by
Extension 2 (applied to $U$ and $T$), $w$ is nontrivial.

In case 2, since $\rho(p)$ and $\rho(q)$ are at least 2 our group
contains rotations by $\pi/2$ about the $x$ and $\ell$ axes.  Since
$\rho(m)\ge 3$, our group also contains $R_z^{\pi/2}$ (as in the proof
of Theorem 4).  Conjugating these rotations about the $x$ and $\ell$
axes by $R_z^{\pi/2}$ gives rotations by $\pi/2$ about the $y$ and
$\ell'$ axes, where $\ell'$ is orthogonal to $\ell$.  Conjugating $A$
and $B$ by these last rotations, we get rotations about the $z$ axis
by $2 \pi/p$ and $2 \pi/q$, and hence by $2 \pi / [p,q]$.  Since
$U^2=R_z^{4\pi/m}$ is also in the group, and since $\rho(m) >
\rho([p,q])$, our group contains $R_z^{2 \pi/r}$, where $r$ is the 
least common multiple $[p,q,m/2]$.  Conjugating this by $R_y^{\pi/2}$
and by $R_{\ell'}^{\pi/2}$ gives $R_x^{2 \pi/r}$ and $R_\ell^{2
\pi/r}$.  Finally, since $n/m$ is an odd multiple of $1/2r$,
conjugating $R_\ell^{2 \pi/r}$ by an appropriate power of $R_z^{2
\pi/r}$ gives $R_{\ell''}^{2 \pi/r}$, where the axis $\ell''$ makes an
angle of $2 \pi/2r$ with the $x$-axis.  So $G_{1/2r}(r,r)\subseteq
G_{n/m}(p,q)$, and therefore $G_{n/m}(p,q)=G_{1/2r}(r,r)$.

Thus $G_{n/m}(p,q)$ does not depend on $p$, $q$, $n$ and $m$
separately, but only on $r$, which is an odd multiple of $m/2$. So we
can replace $p$ and $q$ by 4, $n$ by 1 and $m$ by $2r$, which proves
the first equality of (26).

To prove the second equality, we again assume that $p=q=4$, $n=1$ and
$r=m/2$. Our group is now generated by $A=T$ and $B=U T U^{-1}$, where
$T = R_x^{\pi/2}$ and $U=R_z^{2 \pi/m} = R_z^{\pi/r}$.  Our group also
contains $U^2 = A^2 B^2$.  We consider the subgroups generated by
$U^2$ and $A=T$, on the one hand, and by $U^2$ and $B$ on the other
hand.  Each group is isomorphic to $G(r,4)$, and they have a common
$D_r$ subgroup generated by $U^2$ and $A^2$ (or $U^2$ and $B^2=A^2
U^2$). The group generated by $B$ and $A$ is the quotient of $G(r,4)
*_{D_r} G(r,4)$ by any additional relations.  We will show that there
are no additional relations.

Using (2) and (3) as in the proof of Theorem 1, any element of the
first $G(r,4)$ can be written as
$$ \alpha = W T^{\pm 1} U^{a_1} T^{\pm 1} \cdots T^{\pm 1} U^{a_N}, 
\eqno(30) $$ 
where $W$ is an element of $D_r$, each of the $a_i$'s is even, and,
for $i<N$, $U^{a_i}$ is not a rotation by a multiple of $\pi/2$.
Similarly, an element of the second $G(r,4)$ can be written as
$$ \beta = W U T^{\pm 1} U^{a_1} T^{\pm 1} \cdots T^{\pm 1} U^{a_N-1}. 
\eqno(31) $$ 

We now need a simple fact about amalgamated free products. Consider
some $G *_H \tilde G$, and left coset spaces $G/H$ and $\tilde G/H$.
Pick a set of representatives $\{ k_i \}$ (or $\{ \tilde k_i \}$) for
the left cosets of $G$ (or $\tilde G$).  That is, every element $g$
of $G$ can be written as $g = h k_i$ for some representative $k_i$ and
some $h \in H$, with a similar expression for elements of $\tilde G$.
We claim that any element of $G *_H \tilde G$ can be expressed in the
form:
$$hg_1\tilde g_1g_2\tilde g_2g_3\cdots\tilde g_p\eqno(32)$$
where $h\in H$, each $g_j \in \{k_i\}$ and each $\tilde g_j \in
\{\tilde k_i\}$. To see this, just note that any element of $G *_H
\tilde G$ can automatically be expressed in the form $\phi_1\tilde
\phi_1\phi_2\tilde
\phi_2\phi_3\cdots\tilde \phi_p$, with $\phi_j\in G$ and 
$\tilde \phi_j\in \tilde G$.  Starting from the
right, express $\tilde \phi_p=\tilde h_{p}\tilde k_p$. Then since
$\phi_{p-1}\tilde h_{p} \in G$, we can write $h_{p-1}\tilde
k_{p-1}=\phi_{p-1}\tilde h_{p} $, etc., moving an element of $H$
through the word all the way to the left, where it becomes the $h$ in
(32).

By (30), we may choose representatives for the first $G(r,4)/D_r$ of
the form $T^{\pm 1} U^{a_1} \cdots T^{\pm 1} U^{a_N}$.  Similarly,
(31) indicates the form of the representatives for the second $G(r,4)/D_r$.
So an arbitrary element of $G(r,4) *_{D_r} G(r,4)$ can be written as
$$ \gamma = W (T^{\pm 1} U^{a_{1,1}} T^{\pm 1} \cdots T^{\pm 1} U^{a_{N_1,1}})
U (T^{\pm 1} U^{a_{1,2}} T^{\pm 1} \cdots T^{\pm 1} U^{a_{N_2,2}}) U^{-1}
\cdots, \eqno(33)
$$
where again $W$ is an arbitrary element of $D_r$. Notice that none of
the powers of $U$ are rotations by multiples of $\pi/2$, except
possibly the last term in each parenthetical expression, which then
gets multiplied by $U^{\pm 1}$, yielding an odd power of $U$, hence
not a rotation by a multiple of $\pi/2$. Since $T$ is a rotation by
$\pi/2$, and none of the powers of $U$ are, Extension 2 implies the
resulting word is not equal to the identity, which completes case 2.

In case 3, as in case 2, we have rotations by $\pi/2$ about the $x$
and $z$ axes, hence around the $y$ axis, and hence can transfer
generators from the $x$ to the $z$ axis and vice-versa.  Defining
$U=R_z^{\pi/s}$, as in case 2 we again have $U^2$ in our
group, and, by conjugating by a power of $U^2$, we can exchange the
$\ell$ axis for an axis that makes an angle of $\pi/s$ with the $x$
axis.  This shows that $p$ and $m$ do not contribute separately, but
only through their least common multiple $2s$, and that $n$ does not
matter at all.  This establishes the first equality in (27).

For the second equality we assume $p=4$ and $n=1$, so $m= 2s$ and
$U=R_z^{2 \pi/m}$.  The $G(s,4)$ and $D_q$ subgroups are manifest, as
is their common $D_2$ subgroup.  All that remains is to show that
there are no hidden relations among the elements of $G(s,4) *_{D_2}
D_q$.

The rest of the argument goes just as in case 2.  Find
representatives for all the left cosets of $D_2$ in $G(s,4)$ and in
$D_q$, write out an arbitrary element of the amalgamated free product,
and notice that it isn't trivial.

For case 4, we define $U = R_z^{2 \pi/m}$.  Since $\rho(m)=2$,
$R_z^{\pi/2}$ is an odd power of $U$.  Since $n$ is odd, conjugating
$B$ by an even power of $U$ gives a rotation about the $y$ axis, and
conjugating this by $A^{p/4}$ gives a rotation about the $z$ axis.
Thus generators may be transferred back and forth between the $\ell$ and
$z$ axes, demonstrating that $q$ and $m$ only contribute through their
least common multiple $2t$.  Also, the result is independent of $n$ as usual.
This gives the first equality.

For the second equality, assume $q=2$, so $t=m/2$. The same argument
as in cases 2 and 3 again works. \qed
\bigskip
\nd
{\bf \S 4. Conclusions and Remarks}
\bigskip
The group of rotations of three dimensional Euclidean space is a
fundamental object in mathematics and in mathematical models of
science. Starting with a few simple geometrical facts and a single
algebraic lemma, we have determined the relations in increasingly complicated 
subgroups of the rotation group, 
resulting in a complete classification of the 
``generalized dihedral'' subgroups,
those generated by a pair of rotations of finite order, 
about axes that are themselves separated by an angle which is a rational multiple of $\pi$.

As an application, we consider the role of two-generator subgroups 
of $SO(3)$ in the Banach-Tarski paradox (see [W]).  There one uses a pair
of ``independent'' rotations, that is, rotations which generate a free
subgroup of $SO(3)$, to produce interesting subsets of the unit
sphere. Since in a free group there are NO relations, such rotations
must have infinite order, being rotations by irrational angles.

At first glance such rotations would seem different, indeed
complementary, to the subject of this paper.  However, it is not
difficult to exhibit explicit free 2-generator subgroups of many of
the groups considered here.  For example, let $m$ be any positive
integer other than 1, 2, 4 or 8, let $S=R_y^{\pi/2}$, and let
$T=R_x^{2\pi/m}$.  Then the two rotations $A=STST$ and $B=ST^2ST^2$
are independent.  To see this, write an arbitrary word in $A$ and $B$
in terms of $S$ and $T$.  Since neither $T$ nor $T^2$ is a rotation by
a multiple of $\pi/2$, the result is a word in $S$ and $T$ of the form
(5), and is therefore not the identity. (It is not difficult to check
that the word in $S$ and $T$ is always longer than the original word
in $A$ and $B$.) The group generated by $A$ and $B$ is then a free
2-generator subgroup of $\hat G(m,4)$.

\bigskip
\nd
{\bf Acknowledgements.}\ \ We thank John Conway, J.-P.~Serre and John Tate
for helpful discussions.  
This work, as well as [RS1], is an outgrowth of work with John Conway [CR].
\bigskip
\nd
{\bf  References}
\bigskip \nd
[CR] J.\ Conway and C.\ Radin: Quaquaversal tilings and
rotations, {\it Inventiones math.}, to appear.
[Obtainable from the electronic archive: mp\underbar{
}arc@math.utexas.edu] 
\vs.05 \nd
[RS1] C.\ Radin and L.\ Sadun: Subgroups of
$SO(3)$ associated wtih tilings, {\it J.\ Algebra}, to appear.
[Obtainable from the electronic archive: 
mp\underbar{ }arc@math.utexas.edu]
\vs.05 \nd
[RS2] C.\ Radin and L.\ Sadun: An algebraic invariant for
substitution tiling systems, Geometriae Dedicata, to appear. [Obtainable from the electronic archive:
mp\underbar{ }arc@math.utexas.edu]
\vs.05 \nd
[W] S.\ Wagon: The Banach-Tarski Paradox. Cambridge: The University Press 1985
\vfill\eject
\hbox{}
\vs.1 \hs-.5
\vbox{\epsfig .7\hsize; 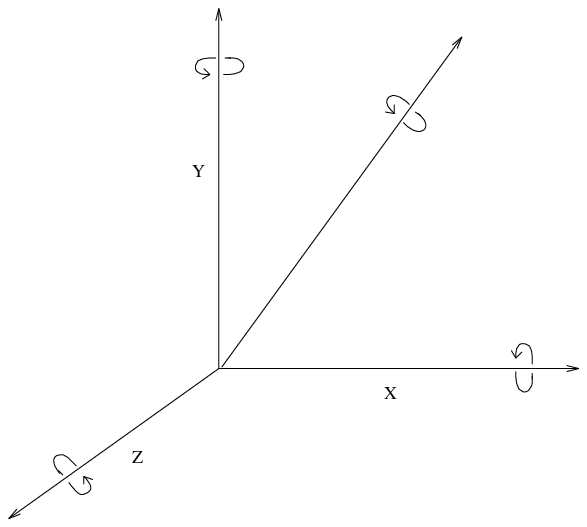;}
\vs-.2 \hs0 $U=R^{2\pi/m}_z$
\vs0 \hs0 $C=R_z^{2\pi/q}$
\vs-2 \hs4.7 $A=R_x^{2\pi/p}$
\vs0 \hs4.7 $T,\hat T=R_x^{\hbox{various angles}}$
\vs-3 \hs3.6 $B=R_\ell^{2\pi/q}$
\vs1 \hs2.7 $\ell$
\vs-1.8 \hs1.7 $S=R_y^{\pi/2}$
\vs5
\centerline{Figure 1}
\vfill \eject

\bye